%%%%%%%%%%%%%%%%%%%%%%%%%%%%%%%%%%%%%%%%%%%%%%%%%%%%%%%%%%%%%%%%%%%%
%AMS-TeX file for the paper
%
%     q-Laguerre polynomials and big $q$-Bessel functions 
%            and their orthogonality relations
%
%by Nicola Ciccoli, Erik Koelink and Tom H. Koornwinder
%%%%%%%%%%%%%%%%%AMS TeXfile%%%%%%%%%%%%%%%%%%%%%%%%%%%%%%%%%%%%%%%%
\input amstex.tex
\documentstyle{amsppt}
\magnification=1200
\baselineskip=13pt
\hsize=6.5truein
\vsize=8.9truein
\parindent=20pt
\NoRunningHeads
\nologo
%%%%%%%%%%%%%%%%% M a c r o s %%%%%%%%%%%%%%%%%%%%%%%%%%%%%%%%%%%%%%
%Section, Theorem and Formula Numbering%%%%%%%%%%%%%%%%%%%%%%%%%%%%%
%%%%%%%%%%%%%%%%%%%%%%%%%%%%%%%%%%%%%%%%%%%%%%%%%%%%%%%%%%%%%%%%%%%%
\countdef\sectionno=1
\countdef\eqnumber=10
\countdef\theoremno=11
\countdef\countrefno=12
\countdef\cntsubsecno=13
\sectionno=0

\def\newsection{\global\advance\sectionno by 1
                \global\eqnumber=1
                \global\theoremno=1
                \global\cntsubsecno=0
                \the\sectionno.\ }

\def\newsubsection{\global\advance\cntsubsecno by 1
                   \the\sectionno.\the\cntsubsecno.\ }

\def\theoremname#1{\the\sectionno.\the\theoremno
                   \xdef#1{{\the\sectionno.\the\theoremno}}
                   \global\advance\theoremno by 1}

\def\eqname#1{\the\sectionno.\the\eqnumber
              \xdef#1{{\the\sectionno.\the\eqnumber}}
              \global\advance\eqnumber by 1}

\global\countrefno=1

\def\refno#1{\xdef#1{{\the\countrefno}}
\global\advance\countrefno by 1}

\def\thmref#1{#1}

%%%%%%%%%%%%%%%%%%%%%%%%%%%%Abbreviations%%%%%%%%%%%%%%%%%%%%%%%%%%%
\def\R{{\Bbb R}}

\def\C{{\Bbb C}}

\def\Z{{\Bbb Z}}
\def\T{{\Bbb T}}

\def\Zp{{\Bbb Z}_{\geq 0}}

\def\hf{{1\over 2}}
\def\ev{{1\over 4}}
\def\al{\alpha}

\def\ga{\gamma}
\def\de{\delta}

\def\et{\eta}
\def\ep{\varepsilon}

\def\la{\lambda}
\def\vp{\varphi}

\def\H{\ell^2(\Z)}
\def\iy{\infty}
\def\FSJ{{\Cal J}}
\def\Pt{\widetilde{P}}
\def\zb{\bar z}
\def\LP{\par\noindent}

\def\LHS{left-hand side}
\def\RHS{right-hand side}

%%%%%%%%%%%%%%%%%%%%%%%%%Reference Numbering%%%%%%%%%%%%%%%%%%%%%%%%
\refno{\Akhi}
\refno{\Bere}
\refno{\Berg}
\refno{\BCG}
\refno{\ChenIM}
\refno{\DunfS}
\refno{\GaspR}
\refno{\Isma}
\refno{\IsmaMS}
\refno{\IsmaR}
\refno{\Kake}
\refno{\KakeMU}
\refno{\KoelITSF}
\refno{\KoelDMJ}
\refno{\KoelIM}
\refno{\KoelV}
\refno{\Koor}
\refno{\KoorS}
\refno{\MassR}
\refno{\MMNNSU}
\refno{\Moak}
\refno{\NM}
\refno{\VK}

%%%%%%%%%%%%%%%%%%%%%%%%%%%%%%%%%%%%%%%%%%%%%%%%%%%%%%%%%%%%%%%%%%%%
%%%%%%%%%%%%%%%%%%%%%%%Beginning of Text%%%%%%%%%%%%%%%%%%%%%%%%%%%%
%%%%%%%%%%%%%%%%%%%%%%%%%%%%%%%%%%%%%%%%%%%%%%%%%%%%%%%%%%%%%%%%%%%%
\topmatter
\title $q$-Laguerre polynomials and big $q$-Bessel functions 
and their orthogonality relations\endtitle 
\author Nicola Ciccoli, Erik Koelink and Tom H. 
Koornwinder\endauthor
%\affil Report W-98-XX, Universiteit van Amsterdam\endaffil
\address Dipartimento di Matematica,
Universit\`a di Perugia, Via Vanvitelli 1, 06123 Perugia, Italy
\endaddress
\email ciccoli\@dipmat.unipg.it\endemail
\address Korteweg-de Vries Instituut, Universiteit van Amsterdam,
Plantage Muidergracht 24, 1018 TV Amsterdam, The Netherlands
\endaddress 
\email koelink\@wins.uva.nl\quad thk\@wins.uva.nl\endemail
\date April 24, 1998\enddate
%\dedicatory dedicated to Dick Askey on his 
%65-th birthday\enddedicatory
\thanks First author
partly supported by C.N.R. - G.N.S.A.G.A.
Second author supported by the Netherlands
Organization for Scientific Research (NWO)
under project number 610.06.100.\endthanks
\abstract The $q$-Laguerre polynomials correspond to an
indetermined moment problem. For explicit discrete non-N-extremal
measures corresponding to Ramanujan's ${}_1\psi_1$-summation
we complement the orthogonal $q$-Laguerre polynomials 
into an explicit orthogonal basis for the corresponding $L^2$-space.
The dual orthogonal system consists
of so-called big $q$-Bessel functions, which can be
obtained as a rigorous limit of the orthogonal system 
of big $q$-Jacobi polynomials. Interpretations on the $SU(1,1)$ 
and $E(2)$ quantum groups are discussed. 
\endabstract 
\endtopmatter
\document

%%%%%%%%%%%%%%%%%%%%%%%%%%%%%%%%%%%%%%%%%%%%%%%%%%%%%%%%%%%%%%%%%%%%
%%%%N E W   S E C T I O N%%%%%%%%%%%%%%%%%%%%%%%%%%%%%%%%%%%%%%%%%%%
%%%%%%%%%%%%%%%%%%%%%%%%%%%%%%%%%%%%%%%%%%%%%%%%%%%%%%%%%%%%%%%%%%%%
%%%%%%%%completely revised introduction%%%%%%%%%%%%%%%%%
\subhead\newsection Introduction\endsubhead

This paper answers two seemingly different questions, which both
originated from quantum groups.
First, consider the system of Moak's \cite{\Moak} 
$q$-Laguerre polynomials
with respect to their familiar discrete orthogonality measure, and
extend this in an explicit way to a complete orthogonal system
of eigenfunctions of a doubly infinite Jacobi matrix originating
from analysis on $SU_q(1,1)$, the quantum $SU(1,1)$ group.
Second, obtain
orthogonality relations and dual orthogonality relations for certain
(big) $q$-Bessel functions originating (see \cite{\BCG}) 
on $E_q(2)$,
the quantum group of plane motions, and give rigorous 
proofs of these
orthogonalities.
The two questions are related
because the dual orthogonality relations for the big 
$q$-Bessel functions
turn out to be the orthogonality relations for the 
completed $q$-Laguerre polynomials.

It is well known that the 
moment problem corresponding to Moak's $q$-Laguerre polynomials
is indetermined as
a Stieltjes moment problem. Moak \cite{\Moak, Thm.~1, 2}
gives several orthogonality measures for the $q$-Laguerre
polynomials: an absolutely continuous measure on $[0,\iy)$
and purely discrete measures supported on the 
set $\{ cq^k\mid k\in\Z\}$ for any constant $c>0$, see \S 4
for the explicit weights. See also Ismail and Rahman 
\cite{\IsmaR} for the explicit calculation 
of the entire functions in the Nevanlinna parametrisation of the
orthogonality measures for the moment problem for the
$q$-Laguerre polynomials.
We are interested in the discrete orthogonality measures. From 
the general theory of orthogonal polynomials it follows 
that the polynomials are not dense in the corresponding
space of quadratically integrable functions, since the 
support is not the set of zeros of an entire function. 
This is expressed
by saying that the measure is not N-extremal, see \cite{\Akhi}
for more information on moment problems and orthogonal polynomials.
In this paper we complement the $q$-Laguerre polynomials to
an orthogonal basis of the $L^2$-space for the discrete measure by
using a certain $q$-analogue of the Bessel function of order
$\al$. The dual basis functions can be recognized as big
$q$-Bessel functions, see below.

$q$-Analogues of Bessel functions exist in several sorts.
Most well-known and probably the oldest ones are Jackson's first and
second $q$-Bessel functions, see Ismail \cite{\Isma}.
They occur in many places, including the present paper, 
but they have
the draw-back that they do not form an orthogonal system 
(possibly they form
a biorthogonal system). Other $q$-analogues of Bessel 
functions can be
obtained as formal limit cases of the three $q$-analogues of
Jacobi polynomials, i.e., of little $q$-Jacobi polynomials,
big $q$-Jacobi polynomials and Askey-Wilson polynomials. For this
reason we propose to speak about little $q$-Bessel functions,
big $q$-Bessel functions and AW type $q$-Bessel functions for
the corresponding limit cases. Little, big and AW type
$q$-Bessel functions have interpretations on $E_q(2)$,
completely analogous to the interpretations of little and big
$q$-Jacobi polynomials and Askey-Wilson polynomials on $SU_q(2)$,
the quantum $SU(2)$ group, 
see Vaksman and Korogodski\u\i\  \cite{\VK}, Koelink 
\cite{\KoelDMJ}, \cite{\KoelIM}, Bonechi et al. 
\cite{\BCG}. See also Ismail et al. \cite{\IsmaMS} for the
Fourier-Bessel transform for the AW type $q$-Bessel function. 
The duals of the AW-type $q$-Bessel functions can also be viewed
as the little $q$-Jacobi functions, which live on $SU_q(1,1)$,
see Masuda et al.  \cite{\MMNNSU}.
These three types
of $q$-Bessel functions satisfy orthogonality relations which
can be obtained as formal limits of the corresponding orthogonality
relations for $q$-Jacobi polynomials. In the $q$-Bessel 
case it is not
sufficient to give orthogonality relations, but one also has to
prove completeness, either directly or by giving dual orthogonality
relations. For rigorous proofs of orthogonality and 
completeness we mention
three different techniques: (i) spectral theoretic methods,
(ii) direct proofs by use of generating functions, 
(iii) rigorous limit
transitions from the orthogonal polynomials case. Usually, only
the first method yields both orthogonality and completeness.
The second method was used for little $q$-Bessel functions by
Koornwinder and Swarttouw \cite{\KoorS}, and was sufficient there
because of self-duality.
The first method was used for little $q$-Jacobi functions by
Kakehi, Masuda and Ueno \cite{\Kake}, \cite{\KakeMU}.
In the present paper we give two different proofs,
first by the spectral method,
and next by method (ii) for the extension of the $q$-Laguerre
orthogonality and by method (iii) for the big 
$q$-Bessel orthogonality.

The spectral method, developed in sections
2  and 3, is based on the spectral analysis of a doubly 
infinite Jacobi matrix, see e.g. Masson and Repka \cite{\MassR}. 
This operator arises in a natural way
from $SU_q(1,1)$. The spectral analysis is
very similar to the spectral analysis for second order 
differential equations  (see \cite{\DunfS}), such as 
for the second order 
differential equation satisfied by the Jacobi functions
(see a survey of these functions in \cite{\Koor}). 
The support of the 
spectral measure is determined by the zeros of a $c$-function. 
We calculate the Green function and we obtain the spectral
measure for the doubly infinite Jacobi matrix from the 
Green function.
The spectral measure is discrete and its support falls  
into two sets; one set corresponding to eigenvectors
in terms of the $q$-Laguerre polynomials and the other
set corresponding to eigenvectors in terms of the $q$-Bessel
coefficients. 
It should be noted that this method corresponds to the one
employed by Kakehi, Masuda and
Ueno \cite{\Kake}, \cite{\KakeMU} giving the Plancherel
measure for the little $q$-Jacobi functions from its
interpretation as matrix elements of irreducible
representations of $SU_q(1,1)$.

The orthogonality
relations and squared norms resulting from the spectral analysis
are explicitly given in \S 4. 
Since we give a basis for this space we also obtain the
dual orthogonality relations as an immediate consequence.
Using a method of Berg \cite{\Berg} we can easily
construct more non-N-extremal orthogonality
measures for the $q$-Laguerre polynomials on the same set
$\{ cq^k\mid k\in\Z\}$ by perturbing with any of the 
$q$-Bessel functions, which is bounded on this set.

In \S 5 we give a straightforward proof of the 
orthogonality relations for the $q$-Laguerre polynomials
and $q$-Bessel functions using two generating functions
for the $q$-Bessel functions. This technique is
motivated by \cite{\KoorS}, \cite{\KoelITSF}. 

In \S 6 we obtain the orthogonality relations for the big $q$-Bessel
functions as a rigorous limit of the big $q$-Jacobi polynomial case.
The formal limit is suggested by the limit transition of
$SU_q(2)$ group to $E_q(2)$.
Big $q$-Jacobi polynomials have an interpretation
as basis elements for the regular representation of $SU_q(2)$
on quantum spheres, see \cite{\NM}.
In the limit transition these basis elements tend 
to the corresponding
basis elements for the regular representation of $E_q(2)$
on quantum hyperboloids, and these latter basis 
elements can be written as
big $q$-Bessel functions, see \cite{\BCG}.

So we see that the same result on $q$-special functions
can be obtained from two different quantum group interpretations:
by considering the quantum $SU(1,1)$ group or the quantum
group of plane motions. This is also the case for the 
little $q$-Jacobi functions studied in 
\cite{\Kake}, \cite{\KakeMU}, which originated as
spherical functions on $SU_q(1,1)$, but which can also be viewed as
certain $q$-analogues of Bessel functions and then 
have an interpretation
on the quantum group of plane motions, see \cite{\KoelIM}.

\demo{Notation} We follow the notation of Gasper and Rahman
\cite{\GaspR} for basic (or $q$-)hypergeometric series.
Throughout we assume that $0<q<1$.
\enddemo

%%%%%%%%%%%%%%%%%%%%%%%%%%%%%%%%%%%%%%%%%%%%%%%%%%%%%%%%%%%%%%%%%%%%
%%%%N E W   S E C T I O N%%%%%%%%%%%%%%%%%%%%%%%%%%%%%%%%%%%%%%%%%%%
%%%%%%%%%%%%%%%%%%%%%%%%%%%%%%%%%%%%%%%%%%%%%%%%%%%%%%%%%%%%%%%%%%%%
\subhead\newsection Solutions to a symmetric operator\endsubhead

Consider the unbounded operator $L$ acting on $\H$ by
$$
\gathered
\bigl( Lu\bigr)_k = a_k\, u_{k+1} + b_k\, u_k +a_{k-1}\, u_{k-1},
\qquad u = (u_k)_{k\in \Z} , \\
a_k = q^{-\hf(k+1)}\sqrt{ 1+c^{-1}q^{-k}}, \qquad
b_k =\sqrt{c^{-1}} (t+t^{-1})q^{-k},
\endgathered
\tag\eqname{\vgldefL}
$$
where $c>0$ and $t\in \R\backslash\{0\}$ are fixed constants,
so that $a_k>0$ and $b_k\in\R$. The operator is densely defined
and symmetric. Split the operator $L$ into two Jacobi matrices $J_+$
and $J_-$, see \cite{\MassR}. 
The coefficients $a_k$ are bounded as $k\to-\iy$,
so that the moment problem corresponding to $J_-$ is determined.
So the deficiency indices for $L$ are either $(0,0)$ or $(1,1)$.
Take $x\in\C\backslash\R$. From the theory of orthogonal 
polynomials we see that the space
of solutions of $Lu=x\, u$ which are
$\ell^2$ for $k\to-\iy$ is one-dimensional.
For the space 
of solutions of $Lu=x\, u$ which are
$\ell^2$ for $k\to\iy$ there are two possibilities:
(i) the space is one-dimensional if the moment
problem for $J_+$ is determined and in that case $L$ is 
self-adjoint;
(ii) the space is two-dimensional if the moment
problem for $J_+$ is indetermined and in that case 
$L$ has a one-parameter family of self-adjoint extensions.
See Akhiezer \cite{\Akhi}, Dunford and Schwartz 
\cite{\DunfS, Ch.~XII} and Masson and Repka \cite{\MassR}
for more information.

Note that, for $k\to\iy$,
$$
a_k\pm b_k+a_{k-1} = c^{-\hf}q^{-k}\bigl(q^\hf+q^{-\hf}
\pm(t+t^{-1})\bigr) + q^{-\hf}\sqrt{c} + {\Cal O}(q^{k/2})
$$
is bounded from above if $\mp t\geq q^{-\hf}$
or $0<\mp t\leq q^\hf$. Hence
we can use the criterion in \cite{\Bere, Ch.~VII} to see
that $J_+$ 
corresponds to a determined moment problem,
and hence $L$ is self-adjoint, if
$|t|\geq q^{-\hf}$ or $0<|t|\leq q^\hf$.

\demo{Remark \theoremname{\remmotivationL}} The motivation
to consider the operator $L$ comes from an investigation 
in the Hopf $\ast$-algebra related to the quantum $SU(1,1)$
group. This is a $\ast$-algebra with two generators $\al$ 
and $\ga$ such that for $0<q<1$ the following relations hold:
$$ 
\al\ga=q\ga\al,\quad 
\al\ga^\ast=q\ga^\ast\al, \quad \ga\ga^\ast=\ga^\ast\ga,
\quad \al\al^\ast-q^2\ga\ga^\ast=1=\al^\ast\al-\ga^\ast\ga.
$$
We can represent this $\ast$-algebra in terms of unbounded
operators in $\H$ by 
$$
\pi(\ga)\, e_k = \la q^k\, e_k, \qquad
\pi(\al)\, e_k = \sqrt{1+|\la|^2q^{2k}}\, e_{k-1},
$$
for the standard orthonormal basis $\{ e_k\}_{k\in\Z}$
of $\H$ and where $\la\in\C\backslash\{0\}$.
Study of the self-adjoint element $\rho=\al^\ast\ga^\ast+\ga\al+
(t+t^{-1})\ga\ga^\ast$ in these representations leads to the
operator in \thetag{\vgldefL} after some normalisation. 
Then we might expect that the
spectral resolution of $L$ will give us information on how
to construct a possible Haar functional on the subalgebra
generated by $\rho$. For the quantum $SU(2)$ group this is
completely rigorous, cf. \cite{\KoelV, \S 5}, and we may 
consider \thetag{\vgldefL} as a non-terminating version of
the three-term recurrence relation for the 
orthonormal Al-Salam and Carlitz polynomials.
\enddemo

For $x\in\C\backslash\{0\}$ define
$$
V^{t}_k(x) = \sqrt{(-q^{1-k}/c;q)_\iy}\;q^{\ev k(k+1)}
\bigl(-t\sqrt{c}\bigr)^k  (qt^2;q)_\infty \ {}_1\vp_1 
\left( {{-t\sqrt{c}/x}\atop{qt^2}};q, xtq^{k+1}\sqrt{c}\right),
\tag\eqname{\eqy}
$$
$$
U_k(x)= \sqrt{(-q^{1-k}/c;q)_\iy}\;q^{\ev k(k+1)}
x^k \
{}_2\vp_1 \left( {{ -\sqrt{c}/(tx), -\sqrt{c}t/x}\atop{0}};
q, -{{q^{1-k}}\over{c}}\right).
\tag\eqname{\eqz}
$$
By \cite{\GaspR, (4.3.2)} we have
for $\pm t\notin\{q^{\hf m}\mid m\in\Z\}$:
$$
\gathered
U_k(x) = C_t\, c_t(x) \, V^t_k (x) + 
C_{t^{-1}}\, c_{t^{-1}}(x) \, V^{t^{-1}}_k (x),\\
c_t(x) = (-\sqrt{c}/xt, qt/x\sqrt{c}, x\sqrt{c}/t;q)_\iy,
\qquad C_t^{-1} = (qt^2,t^{-2}, -c,-q/c;q)_\iy.
\endgathered 
\tag\eqname{\vglconnectieform}
$$
In case $\pm t\in\{q^{\hf m}\mid m\in\Z\}$ we use that for $p\in\Z$
we have
$$
(q^{1-p};q)_\iy\ {}_1\vp_1\left( {{aq^{-p}}\atop{q^{1-p}}};q,z
\right) = {{(q/a;q)_\iy}\over{(q^{p+1}/a;q)_\iy}}
(azq^{-1})^p (q^{1+p};q)_\iy \
{}_1\vp_1\left( {{a}\atop{q^{1+p}}};q,zq^p
\right).
\tag\eqname{\vglestimatektominusinfty}
$$
Formula (\vglestimatektominusinfty) is meaningful and 
can be proved by
interpreting its \LHS\ for  $p\in{\Bbb Z}_{>0}$ as
$$
\sum_{k=0}^\iy{(aq^{-p};q)_k\,(q^{1-p+k};q)_\iy\,
q^{\hf k(k-1)}\,(-z)^k\over (q;q)_k}\,.
$$
{}From \thetag{\vglestimatektominusinfty} we obtain that
$$
V^{\pm q^{-\hf m}}_k(x) = (-c)^m 
{{(\mp q^{1-\hf m}x/\sqrt{c};q)_\infty}
\over {(\mp q^{1+\hf m}x/\sqrt{c};q)_\infty}} 
V^{\pm q^{\hf m}}_k(x).
\tag\eqname{\vgllinafhank}
$$

\proclaim{Lemma \theoremname{\lemsolutions}}
Let $x\in\C\backslash\{0\}$, then $V^t(x)=\bigl( 
V^t_k(x)\bigr)_{k\in\Z}$,
$V^{t^{-1}}(x)=\bigl( V^{t^{-1}}_k(x)\bigr)_{k\in\Z}$ and
$U(x)=\bigl( U_k(x)\bigr)_{k\in\Z}$ are solutions to $Lu = x\, u$. 
Furthermore, 
$V^t(x)$ is $\ell^2$ as $k\to\iy$ if and only if $|t|<q^{-\hf}$
or $t=\pm q^{-\hf m}$, $m\in\Zp$, and
$V^{t^{-1}}(x)$ is $\ell^2$ as $k\to\iy$ if and only if  $|t|>q^\hf$
or $t=\pm q^{\hf m}$, $m\in\Zp$, and $U(x)$ is $\ell^2$ 
as $k\to-\iy$ 
for all $t$. 
If we moreover assume 
$t\in\R\backslash\{q^{\hf m}\mid m\in\Z\}$, then 
$U(x)$ is $\ell^2$ as $k\to\iy$ if and only if 
$q^\hf<|t|<q^{-\hf}$. 
\endproclaim

\demo{Proof} Recall the second order $q$-difference equation
$$
(c-abz)f(qz)+(-(c+q)+(a+b)z)f(z)+(q-z)f(z/q)=0
$$
satisfied by $f(z)={}_2\vp_1(a,b;c;q,z)$, see 
\cite{\GaspR, Exercise 1.13}.
Take $c=0$ to see that $U(x)$ satisfies $Lu=x\,u$. It is clear
from \thetag{\eqz} that $U(x)$ is $\ell^2$ as $k\to -\iy$. 
By confluent limit, $f(z)={}_1\vp_1(a;c;q,z)$ satisfies
$$
(c-az)f(qz)+(-(c+q)+z)f(z)+qf(z/q)=0.
$$
This yields that $V^{t^{\pm 1}}(x)$ satisfies $Lu=x\,u$.
Now use the theta-product identity
$$
(aq^k,q^{1-k}/a;q)_\iy = (-a)^{-k} q^{-\hf k(k-1)} 
(a,q/a;q)_\iy, \qquad \forall\, k\in\Z,\ a\in\C\backslash\{0\},
\tag\eqname{\vglthetaproduct}
$$
to derive the statements on the $\ell^2$-behaviour 
as $k\to\iy$ of $V^{t^{\pm 1}}(x)$ together with 
\thetag{\vgllinafhank}.

Use (\vglconnectieform) to obtain the $\ell^2$-behaviour 
as $k\to\iy$ of $U(x)$.
\qed\enddemo

For $x=0$ we
can take the limit in $V_k^t(x)$, e.g.
$$
V_k^t(0) = \sqrt{(-q^{1-k}/c;q)_\iy}\; q^{\ev k(k+1)}
\bigl(-t\sqrt{c}\bigr)^k  (qt^2;q)_\infty \ {}_0\vp_1 
(-;qt^2;q, -q^{k+1}ct^2),
$$
which is closely related to Jackson's $q$-Bessel function
$J^{(2)}_\al(2\sqrt{c}\,q^{\hf k};q)$ for $t^2=q^\al$, see
Ismail \cite{\Isma}, and also \cite{\ChenIM}, \cite{\Moak}.
Explicitly,
$$
V_k^{q^{\hf\al}}(0) = \sqrt{(-q^{1-k}/c;q)_\iy}\; q^{\ev k(k+1)}
(-1)^k c^{\hf(k-\al)} (q;q)_\iy
J^{(2)}_\al(2\sqrt{c}\,q^{\hf k};q)
\tag\eqname{\vglVKtweenulenJacksonqBessel}
$$
and a similar expression for $V_k^{q^{-\hf\al}}(0)$. We recall the
definition of Jackson's $q$-Bessel function, see \cite{\Isma}:
$$
J^{(2)}_\al(x;q) := {{(q^{\al+1};q)_\iy}\over{(q;q)_\iy}}
\left({{x}\over 2}\right)^\al \ 
{}_0\vp_1\left( {{-}\atop{q^{\al+1}}};q, -q^{\al+1} {{x^2}\over{4}}
\right).
\tag\eqname{\vgldefJacksonqBesselfunction}
$$

\demo{Remark} It follows from Lemma \thmref{\lemsolutions}(i)
that $L$ has deficiency indices $(1,1)$ for $q^\hf<|t|<q^{-\hf}$,
$t\not= \pm 1$, 
and that the deficiency space $N_{\pm i}=
\{ u\in\H\mid Lu=\pm i\, u\}$ is spanned by $U(\pm i)$.
\enddemo

The Wronskian for two sequences $(u_k)_{k\in\Z}$ and 
$(v_k)_{k\in\Z}$ is given by 
$$
[u,v]_k := a_k\bigl( u_{k+1}v_k - u_k v_{k+1}\bigr).
\tag\eqname{\vgldefWronskian}
$$
It is independent of $k$ if $u$ and $v$ are 
two solutions to $Lu =x\, u$, so that we may denote it by $[u,v]$.

\proclaim{Lemma \theoremname{\lemWronskians}} We have the
following Wronskians: 
$$
\align
[V^t(x), V^{t^{-1}}(x)] & = 
{{\sqrt{c}}\over t} (t^2,qt^{-2}-1/c,-cq;q)_\iy, \\
[U(x),V^t(x)] &= 
{{c_{t^{-1}}(x)}\over{-t\sqrt{c}}}.
\endalign
$$
\endproclaim

\demo{Proof} The first Wronskian follows by combining
\thetag{\vgldefWronskian} with Lemma \thmref{\lemsolutions}(ii)
and letting $k\to\iy$ using \thetag{\vglthetaproduct}. 
The second Wronskian follows from the first
and \thetag{\vglconnectieform} for $\pm t\notin 
\{ q^{\hf m}\mid m\in \Z\}$.
Next use analytic continuation with respect to $t$. 
\qed\enddemo

So we find that for $x\in\C\backslash\R$ 
the solution $U(x)$ and $V^{t^{-1}}(x)$  of $Lu =x\, u$ 
are linearly independent for $x\in\C\backslash\R$. 
Moreover, $U(x)$ is $\ell^2$ for $k\to-\infty$ and this determines
$U(x)$ up to a constant by the considerations in the first
paragraph of this section.
Furthermore $V^{t^{-1}}(x)$
is $\ell^2$ for $k\to\infty$ for $|t|>q^\hf$, and for 
$|t|\geq q^{-\hf}$
this condition determines $V^{t^{-1}}(x)$ up to a constant. 

%%%%%%%%%%%%%%%%%%%%%%%%%%%%%%%%%%%%%%%%%%%%%%%%%%%%%%%%%%%%%%%%%%%%
%%%%N E W   S E C T I O N%%%%%%%%%%%%%%%%%%%%%%%%%%%%%%%%%%%%%%%%%%%
%%%%%%%%%%%%%%%%%%%%%%%%%%%%%%%%%%%%%%%%%%%%%%%%%%%%%%%%%%%%%%%%%%%%
\subhead\newsection Spectral resolution\endsubhead

In case $|t|>q^{\hf}$ we see
$V^{t^{-1}}(x)$ is $\ell^2$ for $k\to\iy$ by 
Lemma \thmref{\lemsolutions}(ii). From now on we assume that
$|t|>q^{\hf}$. 
In case $|t|\geq q^{-\hf}$ 
we moreover have that  
$L$ is self-adjoint. 
The domain of $L$ is given by $\{ u\in \H\mid
Lu\in\H\}$. In case, $q^{\hf}<|t|<q^{-\hf}$, $t\not=\pm 1$,
the deficiency indices of $L$
are $(1,1)$, so that all extensions of $L$ to a self-adjoint
operator are parametrised by $U(1)=\T$. In this case, we fix the
self-adjoint extension such that $V^{t^{-1}}$ is contained in the
domain, which is always possible since $\overline{N_i}=N_{-i}$
because of $a_k,b_k\in\R$.
See Dunford and Schwartz \cite{\DunfS, Ch.~XII} for more
information.

Define the Green function for $x\in\C\backslash\R$ by
$$
G_x(m,n) := \cases {\displaystyle{
{{U_n(x) V^{t^{-1}}_m(x)}\over{[U(x), V^{t^{-1}}(x)]}},}}
&\text{for $n\leq m$,}\\
{\displaystyle{{{U_m(x) V^{t^{-1}}_n(x)}\over{[U(x), 
V^{t^{-1}}(x)]}},}}
&\text{for $n> m$,} \endcases
$$
so that $\bigl( (x-L)^{-1} f\bigr)_m = \sum_{n=-\iy}^\iy
G_x(m,n)\, f_n$. For the resolution of the identity
$E$ for $L$, i.e. $L= \int_\R t\, dE(t)$, we now have,
see \cite{\DunfS, Ch.~XII},
$$
\langle E\bigl( (a,b)\bigr) v,w\rangle 
= \lim_{\de\downarrow 0}\lim_{\ep\downarrow 0}
{1\over{2\pi i}} \int_{a+\de}^{b-\de} 
\langle (s-i\ep-L)^{-1}v,w\rangle - 
\langle (s+i\ep-L)^{-1}v,w\rangle \, ds.
$$

Now observe that
$$
\langle (s\pm i\ep-L)^{-1}v,w\rangle =
\sum_{n\leq m} {{U_n(s\pm i\ep)V^{t^{-1}}_m(s\pm i\ep)}\over{
[U(s\pm i\ep),V^{t^{-1}}(s\pm i\ep)]}} {{v_n \bar w_m + v_m\bar w_n}
\over 2}.
$$
Since $V^{t^{-1}}_m(x)$ is an entire function of $x$, and
$U_n(x)$ is analytic in $\C\backslash\{ 0\}$ we see that
the measure $\langle E\bigl( (a,b)\bigr) v,w\rangle$ has only
discrete mass points at the zeros of the Wronskian
$[U(x), V^{t^{-1}}(x)]$. By Lemma \thmref{\lemWronskians}
this corresponds to the zeros of $c_t(x)$ which are the
points $\et_p = -\sqrt{c}q^p/t$, $p\in\Zp$, and
$\xi_p=tq^p/\sqrt{c}$, $p\in\Z$. So the spectrum
of $L$ is $\{\et_p\mid p\in \Zp\} \cup \{\xi_p\mid p\in \Z\}
\cup \{ 0 \}$, and it remains to show that
$0$ is not in the point spectrum 
of $L$. It follows from Moak \cite{\Moak, \S 7} and 
\thetag{\vglVKtweenulenJacksonqBessel} that $V_k^{t^{-1}}(0)$
is not identically zero for $k<-N$ for some $N$. The
observation now follows from the asymptotic expression
for the Jackson $q$-Bessel function derived by
Chen et al. \cite{\ChenIM, Thm.~4}, which shows that the
Jackson $q$-Bessel function in 
\thetag{\vglVKtweenulenJacksonqBessel} increases exponentially
with $k^2$ as $k\to -\iy$. 

In order to calculate $E(\{\et_p\})$ we take 
the interval $(a,b)$ such that it contains only $\et_p$
as a point from the spectrum. Then
$$
\align
\langle E\bigl( \{\et_p\}\bigr) &v,w\rangle =
\langle E\bigl( (a,b)\bigr) v,w\rangle \\ &= 
\lim_{\ep\downarrow 0}
{1\over{2\pi i}} \int_a^b
\langle (s-i\ep-L)^{-1}v,w\rangle - 
\langle (s+i\ep-L)^{-1}v,w\rangle \, ds \\ &=
{1\over {2\pi i}} \oint_{(\et_p)} 
\langle (s-L)^{-1} v,w\rangle\, ds  \\ &=
\sum_{n\leq m}{{v_n \bar w_m + v_m\bar w_n}
\over 2} {1\over {2\pi i}} \oint_{(\et_p)}  
{{U_n(s)V^{t^{-1}}_m(s)}\over{
[U(s),V^{t^{-1}}(s)]}} \, ds\\ &=
\sum_{n\leq m}{{v_n \bar w_m + v_m\bar w_n}
\over 2} U_n(\et_p) V^{t^{-1}}_m(\et_p) 
(-\sqrt{c}t^{-1}) 
\text{Res}_{x=\et_p} {1\over{c_t(x)}} \\ &=
C_{t^{-1}} c_{t^{-1}}(\et_p) (-\sqrt{c}t^{-1}) 
\text{Res}_{x=\et_p} {1\over{c_t(x)}}
\langle v, V^{t^{-1}}(\et_p)\rangle \langle  
V^{t^{-1}}(\et_p),w\rangle 
\endalign
$$
by \thetag{\vglconnectieform} and $c_t(\et_p)=0$. 
Since ${}_1\vp_1(q^{-p};qt^{-2};q, -cq^{p+k+1}t^{-2})
\sim C q^{kp}$ for $k\to-\iy$, we obtain $V^{t^{-1}}(\et_p)\in\H$
also in a direct manner. Next 
$$
\text{Res}_{x=\et_p} {1\over{c_t(x)}} = 
{1\over{(-q^{1-p}t^2/c,-q^pc/t^2;q)_\iy}}
{{-q^p\sqrt{c}/t}\over{(q^{-p};q)_p(q;q)_\iy}}
$$
and
$$
C_{t^{-1}} c_{t^{-1}}(\et_p) = 
{{(t^2q^{-p};q)_p }\over{(qt^{-2};q)_\infty}}
c^{-p} q^{-\hf p(p-1)}
$$
by \thetag{\vglthetaproduct}. 
So that finally,
$$
\langle E\bigl( \{\et_p\}\bigr) v,w\rangle = 
c t^{-2} q^p {{(qt^{-2};q)_p}\over{(q;q)_p}}
{1\over{(-c/t^2,-qt^2/c,qt^{-2};q)_\iy}}
\langle v, V^{t^{-1}}(\et_p)\rangle \langle  
V^{t^{-1}}(\et_p),w\rangle.
\tag\eqname{\vglspectralonZplus}
$$
Take $v=w=V^{t^{-1}}(\et_p)$ and use that 
$E\bigl( \{\et_p\}\bigr)V^{t^{-1}}(\et_p) = 
V^{t^{-1}}(\et_p)$ to see 
that 
$$
\| V^{t^{-1}}(\et_p)\|^2 = 
c^{-1} t^2 q^{-p} {{(q;q)_p}\over{(qt^{-2};q)_p}}
(q,-c/t^2,-qt^2/c,qt^{-2};q)_\iy
\tag\eqname{\vglsquarednormonZPlus}
$$
Since the projections $E\bigl( \{\et_p\}\bigr)$ satisfy 
$E\bigl( \{\et_p\}\bigr)E\bigl( \{\et_r\}\bigr)=\de_{pr}
E\bigl( \{\et_p\}\bigr)$ we get
$$
\multline
\langle V^{t^{-1}}(\et_p),V^{t^{-1}}(\et_r)\rangle = 
\langle E\bigl( \{\et_p\}\bigr)V^{t^{-1}}(\et_p),
E\bigl( \{\et_r\}\bigr)V^{t^{-1}}(\et_r) \rangle = \\
\langle V^{t^{-1}}(\et_p), E\bigl( \{\et_p\}\bigr)
E\bigl( \{\et_r\}\bigr)V^{t^{-1}}(\et_r) \rangle =
\de_{pr}\langle V^{t^{-1}}(\et_p),V^{t^{-1}}(\et_p) \rangle
\endmultline
\tag\eqname{\vglorthorelonZplus}
$$
and after setting $t^{-2}=q^\al$, $\al>-1$, we obtain 
from \thetag{\vglsquarednormonZPlus} and 
\thetag{\vglorthorelonZplus} the orthogonality relations
for the $q$-Laguerre polynomials on the set $\{ cq^k\mid k\in\Z\}$,
$c>0$, obtained by Moak \cite{\Moak}. It's well known that these
polynomials correspond to an indetermined moment problem and that
the orthogonality measure supported on $\{ cq^k\mid k\in\Z\}$
is not N-extremal, 
i.e. the polynomials are not dense in the
corresponding $L^2$-space, see \cite{\Akhi}, \cite{\Moak}. 

Similarly,
$$
\langle E\bigl( \{\xi_p\}\bigr) v,w\rangle =
C_{t^{-1}} c_{t^{-1}}(\xi_p) (-\sqrt{c}t^{-1}) 
\text{Res}_{x=\xi_p} {1\over{c_t(x)}}
\langle v, V^{t^{-1}}(\xi_p)\rangle \langle  
V^{t^{-1}}(\xi_p),w\rangle 
$$
with 
$$
\gather
\text{Res}_{x=\xi_p} {1\over{c_t(x)}} = {{tq^p}\over{\sqrt{c}}}
{{(-1)^{p-1}q^{\hf p(p-1)}}\over{(-cq^{-p}t^{-2},q,q;q)_\iy}},\\
C_{t^{-1}} c_{t^{-1}}(\xi_p) = \bigl( - {c\over{t^2}}\bigr)^p
q^{-p^2} {1\over{(-q^{1+p}/c;q)_\iy}}.
\endgather
$$
Hence, 
$$
\langle E\bigl( \{\xi_p\}\bigr) v,w\rangle =
q^p {{(-q^{p+1}t^2/c;q)_\iy}\over{
(-q^{p+1}/c,q,q,-ct^{-2},-qt^2/c;q)_\iy}}
\langle v, V^{t^{-1}}(\xi_p)\rangle \langle  
V^{t^{-1}}(\xi_p),w\rangle .
$$
Then, similarly as before, we obtain
$$
\langle V^{t^{-1}}(\xi_p),V^{t^{-1}}(\xi_r)\rangle = \de_{pr}
q^{-p} 
{{(-q^{p+1}/c,q,q,-ct^{-2},-qt^2/c;q)_\iy}
\over{(-q^{p+1}t^2/c;q)_\iy}}.
\tag\eqname{\vglorthorelonwholeZ}
$$

\demo{Remark}
The vector $V^{t^{-1}}(\xi_p)$ belongs to $\H$ for $p\in\Z$.
This can also be seen as follows. First apply to 
(\eqz) the transformation
formula
$$
{}_1\vp_1\left({a\atop c};q,z\right)=
{(z;q)_\iy\over(c;q)_\iy}\
{}_1\vp_1\left({az/c\atop z};q,c\right)
\tag\eqname{\eqaa}
$$
(a limit case of Heine's
transformation formula \cite{\GaspR, (1.4.5)}).
Combination of (\eqz), (\eqaa) and (
\vglestimatektominusinfty) yields that
$$
\multline
V_k^{t^{-1}}(\xi_p)=
{(-q^{p+1}/c;q)_\iy\,(-t^2c^{-1}q^p)^p\over
\sqrt{(-q^{1-k}/c;q)_\infty} (qt^{-2};q)_\iy}\,
q^{\ev k(k+1)}\,(tc^{-\hf}q^p)^k (qt^{-2};q)_\infty\\ \times\,
(q^{1-p-k};q)_\iy 
\,{}_1\vp_1\left( {{-cq^{-p}}\atop{q^{1-p-k}}};q,q^{1-p-k}
t^{-2}\right),
\endmultline
$$
so that we obtain the $\ell^2$-behaviour as $k\to -\iy$
of $V^{t^{-1}}(x)$.
\enddemo

%%%%%%%%%%%%%%%%%%%%%%%%%%%%%%%%%%%%%%%%%%%%%%%%%%%%%%%%%%%%%%%%%%%%
%%%%N E W   S E C T I O N%%%%%%%%%%%%%%%%%%%%%%%%%%%%%%%%%%%%%%%%%%%
%%%%%%%%%%%%%%%%%%%%%%%%%%%%%%%%%%%%%%%%%%%%%%%%%%%%%%%%%%%%%%%%%%%%
\subhead\newsection Basis for $L^2$-space\endsubhead

By $L^2(\mu^{(\al;c)})$ we denote the space of 
square integrable functions
on the set $\{ cq^k\mid k\in\Z\}$, $c>0$, with positive weight
$q^{k(\al+1)}/(-cq^k;q)_\iy$ at $cq^k$, $k\in \Z$, i.e.
$f\in L^2(\mu^{(\al;c)})$ if 
$$
{\Cal L}(|f|^2) := 
\sum_{k=-\iy}^\iy {{q^{k(\al+1)}}\over{(-cq^k;q)_\iy}}
|f(cq^k)|^2<\iy.
$$
Here we take $\al>-1$. 

Recall the $q$-Laguerre polynomials introduced by Moak \cite{\Moak},
see also \cite{\GaspR},
$$
L_n^{(\al)}(x;q) := {{(q^{\al+1};q)_n}\over{(q;q)_n}}
{}_1\vp_1\left( {{q^{-n}}\atop{q^{\al+1}}};q,-xq^{n+\al+1}\right).
\tag\eqname{\vgldefqLaguerrepol}
$$
We also define the functions
$$
\aligned
M_p^{(\al;c)}(x;q) :=& {{(q^{\al+1};q)_\iy}\over{
(q,-cq^{\al+1};q)_\iy}} \, {}_1\vp_1
\left( {{-cq^{\al-p}}\atop{q^{\al+1}}};q, {{xq^{p+1}}\over{c}}
\right) \\
=& {{(xq^{p+1}/c;q)_\iy}\over{
(q,-cq^{\al+1};q)_\iy}} \, {}_1\vp_1 \left(
{{-x}\atop{xq^{p+1}/c}};q,q^{\al+1}\right)
\endaligned
\tag\eqname{\vgldefrestofbasis}
$$
for $p\in\Z$. (The second equality is by (\eqaa).)

\proclaim{Theorem \theoremname{\thmbasisappropriateLtwospace}}
The functions $M_p^{(\al;c)}(\cdot;q)$, $p\in\Z$, 
together with the $q$-Laguerre polynomials
$L_n^{(\al)}(\cdot;q)$, $n\in\Zp$, form an orthogonal basis
for $L^2(\mu^{(\al;c)})$, $c>0$, $\al>-1$. Explicitly,
$$
\align
{\Cal L}\bigl( L_n^{(\al)}(\cdot;q)L_p^{(\al)}(\cdot;q)\bigr)
&= \de_{n,p} q^{-p} {{(q^{\al+1};q)_p}\over{(q;q)_p}}
{{(q,-cq^{\al+1}, -q^{-\al}/c;q)_\iy}\over
{(q^{\al+1},-c,-q/c;q)_\iy}}, \\
{\Cal L}\bigl( M_p^{(\al;c)}(\cdot;q)M_r^{(\al;c)}(\cdot;q)\bigr)
&= \de_{p,r} cq^\al q^{-p} {{(-q^{p+1}/c,-q^{-\al}/c;q)_\iy}\over
{(-q^{p+1-\al}/c,-cq^{\al+1};q)_\iy}}
 {1\over{(-c,-q/c;q)_\iy}}, \\
{\Cal L}\bigl( M_p^{(\al;c)}(\cdot;q)L_n^{(\al)}(\cdot;q)\bigr)
&= 0.
\endalign
$$
\endproclaim

Note that $M_p^{(\al;c)}(x;q)$ depends on $c$,
unlike $L_n^{(\al)}(x;q)$. It should also be observed
that $M_p^{(\al;c)}(cq^k;q)$ is bounded as $k\to\iy$
and that $M_p^{(\al;c)}(cq^k;q)\to 0$ for $|q^{p-\al}/c|<1$
for $k\to-\iy$ as follows from
\thetag{\vglestimatektominusinfty}. Using the method of
Berg \cite{\Berg} we can construct more orthogonality 
measures for the $q$-Laguerre polynomials supported on
the set $\{ cq^k\mid k\in\Z\}$. Let
then $|M_p^{(\al;c)}(cq^k;q)|\leq K$ for some $K$
for any $p\in\Z$ satisfying $|q^{p-\al}/c|<1$. 
Then the $q$-Laguerre
polynomials are also orthogonal with respect to the 
positive discrete measure with masses 
$\bigl(1+sK^{-1}M_p^{(\al;c)}(cq^k;q)\bigr)
q^{k(\al+1)}/(-cq^k;q)_\iy$ 
at $cq^k$, $k\in \Z$, for any $s\in[-1,1]$, cf.
\cite{\Berg, Prop.~4.1}. We may also take suitable linear
combinations of the functions $M_p^{(\al;c)}(x;q)$.

\demo{Proof} Apply the results of the previous section with
$t^{-2} = q^\al$. The first statement corresponds
to $\langle V^{t^{-1}}(\et_p),V^{t^{-1}}(\et_r)\rangle = \de_{p,r}
\| V^{t^{-1}}(\et_p)\|^2$, the second statement corresponds to
$\langle V^{t^{-1}}(\xi_p),V^{t^{-1}}(\xi_r)\rangle = \de_{p,r}
\| V^{t^{-1}}(\xi_p)\|^2$ and the last one to
$\langle V^{t^{-1}}(\xi_p),V^{t^{-1}}(\et_r)\rangle = 0$.
\qed\enddemo

\demo{Remark \theoremname{\remthmbasisappropriateLtwospace}} 
The first statement of 
Theorem \thmref{\thmbasisappropriateLtwospace} corresponds
to Moak's discrete orthogonality relations \cite{\Moak, Thm.~2}
for the $q$-Laguerre polynomials. The second relation can 
be rewritten as
$$
\multline
\sum_{k=-\iy}^\iy {{q^{k(\al+1)}(q^{\al+1};q)_\iy}
\over{(-cq^k;q)_\iy(q;q)_\iy}}
{}_1\vp_1\left( {{-cq^{\al-p}}\atop{q^{\al+1}}};q,q^{p+k+1}\right) 
{{(q^{\al+1};q)_\iy}\over{(q;q)_\iy}}
{}_1\vp_1\left( {{-cq^{\al-r}}\atop{q^{\al+1}}};q,q^{r+k+1}\right)\\
=  \de_{p,r} cq^\al
q^{-p} {{(-q^{p+1}/c;q)_\iy}\over
{(-q^{p+1-\al}/c;q)_\iy}} 
{{(-cq^{\al+1},-q^{-\al}/c;q)_\iy}\over{(-c,-q/c;q)_\iy}} = \\
\sum_{k=-\iy}^\iy {{q^{k(\al+1)}(q^{k+p+1};q)_\iy}
\over{(-cq^k;q)_\iy (q;q)_\iy}}
{}_1\vp_1\left( {{-cq^k}\atop{q^{k+p+1}}};q,q^{\al+1}\right)
{{(q^{k+r+1};q)_\iy}\over{(q;q)_\iy}}
{}_1\vp_1\left( {{-cq^k}\atop{q^{k+r+1}}};q,q^{\al+1}\right),
\endmultline
\tag\eqname{\vglorthorelationscomplete}
$$
for $p,r\in\Z$,
which can be viewed as a $q$-analogue of the Hankel 
transform 
for the first equality or as a $q$-analogue of the Hansen-Lommel 
orthogonality relations for the second equality,
cf. \cite{\KoorS}, and the limit case $c\to 0$ corresponds
to \cite{\KoorS, Prop.~2.6}.  
The last statement shows that these $q$-Bessel functions are 
orthogonal to the $q$-Laguerre polynomials;
$$
\sum_{k=-\iy}^\iy {{q^{k(\al+1)}}\over{(-cq^k;q)_\iy}}
L^{(\al)}_n(cq^k;q) {{(q^{\al+1};q)_\iy}\over{(q;q)_\iy}}
{}_1\vp_1\left( 
{{-cq^{\al-r}}\atop{q^{\al+1}}};q,q^{r+k+1}\right) =0,
\quad r\in\Z,\, n\in\Zp.
$$
\enddemo

Since the statement is that $\{V^{t^{-1}}(\et_p)\}_{p\in\Zp}$
and $\{V^{t^{-1}}(\xi_p)\}_{p\in\Z}$ form an orthogonal basis
for the Hilbert space $\H$, we also find the dual orthogonality
relations 
$$
\sum_{p=0}^\iy {{V_k^{t^{-1}}(\et_p)V_l^{t^{-1}}(\et_p)}\over
{\| V^{t^{-1}}(\et_p)\|^2}} 
+ \sum_{p=-\iy}^\iy {{V_k^{t^{-1}}(\xi_p)V_l^{t^{-1}}(\xi_p)}\over
{\| V^{t^{-1}}(\xi_p)\|^2}} =\de_{k,l}.
\tag\eqname{\vglformdualorthorelations}
$$
The first sum is the Poisson kernel for the $q$-Laguerre
polynomials evaluated at one; it can also be derived from
the Christoffel-Darboux formula and the limit transition
of the $q$-Laguerre polynomials to Jackson's $q$-Bessel function,
see \cite{\Moak}. Explicitly, from \cite{\Moak, (4.11), 
Thm.~5} we get
$$
\multline
\sum_{p=0}^N {{q^p(q;q)_p}\over{(q^{\al+1};q)_p}}
L^{(\al)}_p(x;q)L^{(\al)}_p(y;q) = \\
{{(q;q)_N}\over{(q^{\al+1};q)_N}} {1\over{x-y}}
\bigl( x L^{(\al+1)}_N(x;q)L^{(\al)}_N(y;q) - 
y L^{(\al+1)}_N(y;q)L^{(\al)}_N(x;q)\bigr) 
{\overset{N\to\iy}\to{\longrightarrow}} \\
{{(q;q)_\iy}\over{(q^{\al+1};q)_\iy}} 
{{(xy)^{-\hf\al}}\over{x-y}}
\bigl( \sqrt{x}
J^{(2)}_{\al+1}(2\sqrt{x};q)J^{(2)}_\al(2\sqrt{y};q) - 
\sqrt{y} J^{(2)}_{\al+1}(2\sqrt{y};q)J^{(2)}_\al(2\sqrt{x};q)\bigr),
\endmultline
\tag\eqname{\eqm}
$$
where the right hand side is well-defined for $x=y$
using l'H\^opital's formula. Using $t^{-2}=q^\al$ we get
$$
\multline
\sum_{p=0}^\iy {{V_k^{t^{-1}}(\et_p)V_l^{t^{-1}}(\et_p)}\over
{\| V^{t^{-1}}(\et_p)\|^2}} =  
{{q^{\hf(\al+1)(k+l)} (-1)^{k+l}}\over
{(-cq^k,-cq^l;q)_\iy^\hf}} {{(-c,-q/c,q^{\al+1};q)_\iy}\over{
(-cq^{\al+1},-q^{-\al}/c,q;q)_\iy}} \\ \times
\sum_{p=0}^\iy {{q^p(q;q)_p}\over{(q^{\al+1};q)_p}}
\, L_p^{(\al)}(cq^k;q)L_p^{(\al)}(cq^l;q) 
= {{q^{\hf(k+l)}(-1)^{k+l}c^{-\al-\hf}(-c,-q/c;q)_\iy}\over
{(-cq^k,-cq^l;q)_\iy^\hf (-cq^{\al+1},-q^{-\al}/c;q)_\iy}}
\\ \times {1\over{q^k-q^l}} \Bigl(
q^{\hf k}J^{(2)}_{\al+1}(2\sqrt{c}q^{\hf k};q)
J^{(2)}_\al(2\sqrt{c}q^{\hf l};q) - 
q^{\hf l}J^{(2)}_{\al+1}(2\sqrt{c}q^{\hf l};q)
J^{(2)}_\al(2\sqrt{c}q^{\hf k};q)\Bigr)
\endmultline
$$
by \thetag{\eqm}.
Calculating the second sum in \thetag{\vglformdualorthorelations}
is straightforward, so that we now have obtained the following
corollary.

\proclaim{Corollary \theoremname{\eql}}
The following orthogonality relations hold;
$$
\multline
\de_{k,l} cq^{-k} (-cq^k;q)_\iy 
{{(-cq^{\al+1},-q^{-\al}/c;q)_\iy}\over{(-c,-q/c;q)_\iy}}
= \\ {{c^{\hf-\al}}\over{q^k-q^l}}
\Bigl(
q^{\hf k}J^{(2)}_{\al+1}(2\sqrt{c}q^{\hf k};q)
J^{(2)}_\al(2\sqrt{c}q^{\hf l};q) - 
q^{\hf l}J^{(2)}_{\al+1}(2\sqrt{c}q^{\hf l};q)
J^{(2)}_\al(2\sqrt{c}q^{\hf k};q)\Bigr)
+ \\ q^{\al(\hf(k+l)-1)} \sum_{p=-\iy}^\iy
q^p {{(-q^{p+1-\al}/c;q)_\iy}\over{(-q^{p+1}/c;q)_\iy}}
{{(q^{\al+1};q)_\iy}\over{(q;q)_\iy}}
{}_1\vp_1\left( {{-cq^{\al-p}}\atop{q^{\al+1}}};q,q^{k+p+1}\right)
\\ \times {{(q^{\al+1};q)_\iy}\over{(q;q)_\iy}}
{}_1\vp_1\left( {{-cq^{\al-p}}\atop{q^{\al+1}}};q,q^{l+p+1}\right),
\endmultline
$$
\endproclaim

In this form Corollary \thmref{\eql}
is reminiscent of the Hankel transform, whereas if we use the
transformation for the ${}_1\vp_1$-series of 
\thetag{\vgldefrestofbasis}, the orthogonality relations
remind us of the Hansen-Lommel orthogonality relations,
cf. \cite{\KoorS}. Of course, we may also replace the 
expression for the Jackson $q$-Bessel function by the
Poisson kernel for the $q$-Laguerre polynomials evaluated
at one, cf. \thetag{\eqm}, to
obtain the dual orthogonality relations involving a sum over
$\Z$ and over $\Zp$.

%%%%%%%%%%%%%%%%%%%%%%%%%%%%%%%%%%%%%%%%%%%%%%%%%%%%%%%%%%%%%%%%%%%%
%%%%N E W   S E C T I O N%%%%%%%%%%%%%%%%%%%%%%%%%%%%%%%%%%%%%%%%%%%
%%%%%%%%%%%%%%%%%%%%%%%%%%%%%%%%%%%%%%%%%%%%%%%%%%%%%%%%%%%%%%%%%%%%
\subhead\newsection Direct proof of the orthogonality 
relations of
Theorem \thmref{\thmbasisappropriateLtwospace}\endsubhead

In this section we present a direct analytic 
proof using summation and transformation formulas
of the orthogonality
relations of Theorem \thmref{\thmbasisappropriateLtwospace}, 
but not of the completeness. The first statement of the
orthogonality relations is well known, see Moak \cite{\Moak},
so we concentrate on the last two. 
The method is based on manipulating generating functions,
cf. \cite{\KoorS}, \cite{\KoelITSF}.
We start with the
following lemma giving generating functions for the 
$q$-Bessel functions under consideration.

\proclaim{Lemma \theoremname{\lemgeneratingfunctions}}
We have the following generating functions: \par\noindent
{\rm (i)} For $0<|z|<|b|^{-1}$ we have
$$
\sum_{p=-\iy}^\iy z^pb^p {{(q^{p+1};q)_\iy}\over
{(aq^p/b;q)_\iy}} \, {}_1\vp_1\left( {{aq^p/b}\atop{q^{p+1}}};
q,bx\right) = {{(q,az,x/z;q)_\iy}\over{(a/b,bz;q)_\iy}},
$$
\par\noindent
{\rm (ii)} For $|d|<|w|<1$ we have
$$
\sum_{r=-\iy}^\iy w^r (q^{r+1};q)_\iy
\, {}_1\vp_1\left( {{dq^{r+1}/y}\atop{q^{r+1}}};
q,y\right) = {{(d,q,y/w;q)_\iy}\over{(w,d/w;q)_\iy}}.
$$
\endproclaim

\demo{Proof} Case (i) is a limit case of 
\cite{\KoelITSF, Prop.~2.1}. Use the $q$-binomial
theorem \cite{\GaspR, (1.3.2)} to expand 
$(az;q)_\iy/(bz;q)_\iy$ in a power series of
$z$ for $|z|<|b|^{-1}$, and \cite{\GaspR, (1.3.16)}
to expand $(x/z;q)_\iy$ in a power series in
$z^{-1}$ for $|z|>0$. Multiply the resulting series
to find the result.

For the proof of (ii) we expand $1/(w;q)_\iy$ in a power 
series in $w$ for $|w|<1$ and 
$(y/w;q)_\iy/(d/w;q)_\iy$ in a power series
in $w^{-1}$ for $|d|<|w|$ using \cite{\GaspR, (1.3.2)}
twice. Combine this to write the right hand side 
as a Laurent series in $w$ with coefficients given by
a ${}_2\vp_1$-series. Then use the limit case $b\to 0$ of Heine's
transformation \cite{\GaspR, (1.4.6)} to get the result.
\qed
\enddemo

To prove the last part of 
Theorem \thmref{\thmbasisappropriateLtwospace}
we use Lemma \thmref{\lemgeneratingfunctions}(i)
after replacing $p$ by $k+r$, $a$ by $-bcq^{-r}$ and
$x$ by $q^{\al+1}/b$ to get
$$
\sum_{k=-\iy}^\iy {{(zb)^{k+r}}\over{(-cq^k;q)_\iy}}
(q^{k+r+1};q)_\iy \, {}_1\vp_1\left( {{-cq^k}\atop{q^{k+r+1}}};
q,q^{\al+1}\right) =
{{(q,-bzcq^{-r},q^{\al+1}/(bz);q)_\iy}\over
{(-cq^{-r},bz;q)_\iy}}
$$
and specialising $zb=q^{\al+1+m}$, $\al>-1$, shows that the right
hand side is zero for $m\in\Zp$ whereas the left hand side
is a non-zero multiple of ${\Cal L}$ applied to
$M^{(\al;c)}_r(x;q)x^m$. So $M^{(\al;c)}_r(x;q)$ is orthogonal
to all monomials, hence to all polynomials implying the
last statement of
Theorem  \thmref\thmref{\thmbasisappropriateLtwospace}. 

To prove the orthogonality relations for the functions
$M^{(\al;c)}(x;q)$ of 
Theorem \thmref{\thmbasisappropriateLtwospace} we first
deduce the following result.

\proclaim{Proposition \theoremname{\propfirststeporthorels}}
For $|y|<1$, $l\in\Z$, we have
$$
\multline
\sum_{k=-\iy}^\iy {{y^k}\over{(aq^k/b;q)_\iy}}
(q^{k+1};q)_\iy {}_1\vp_1\left( {{aq^k/b}\atop{q^{k+1}}};q,y
\right)
(q^{k-l+1};q)_\iy {}_1\vp_1\left( {{dq^{k-l+1}/y}
\atop{q^{k-l+1}}};q,y \right) \\
= \cases 0, &\text{for $l<0$,}\\
{\displaystyle{ d^l{{(ay/(bd);q)_l (d,q,q;q)_\iy}\over{
(q;q)_l(a/b;q)_\iy}}}}, &\text{for $l\geq 0$.}
\endcases
\endmultline
$$
\endproclaim

\demo{Proof} Choose $w=y/(bz)$ in 
Lemma \thmref{\lemgeneratingfunctions}(ii) and multiply
with the generating function of
Lemma \thmref{\lemgeneratingfunctions}(i) with $x$
replaced by $y/b$. This then gives the following
identity
$$
\multline
{{(az;q)_\iy (d,q,q;q)_\iy}\over{
(bzd/y;q)_\iy (a/b;q)_\iy}} =  
\sum_{l=-\iy}^\iy z^l \\ \times
\sum_{k=-\iy}^\iy
{{y^{k-l}b^l}\over {(aq^k/b;q)_\iy}}
(q^{k+1};q)_\iy {}_1\vp_1\left( {{aq^k/b}\atop{q^{k+1}}};q,y
\right)
(q^{k-l+1};q)_\iy {}_1\vp_1\left( {{dq^{k-l+1}/y}
\atop{q^{k-l+1}}};q,y \right)
\endmultline
$$
under the condition $|y/b|<|z|<\min( |b|^{-1},|y/(db)|)$.
So we need $|y|<1$, $|d|<1$ to have a non-empty region of
analyticity of the right hand side as function of $z$.
The left hand side is an analytic function of $z$ in
$|z|<|y/(db)|$ and it can be expanded in a power series in
$z$ using the $q$-binomial theorem \cite{\GaspR, (1.3.2)}.
Comparing coefficients at both sides gives the result 
for $|d|<1$.

Next we use analytic continuation with respect to $d$. For fixed
$l\in\Z$ the right hand side is analytic in $d$. To see that the
left hand side is analytic in $d$ we note that 
$f_k(d)=(q^{k-l+1};q)_\iy {}_1\vp_1(dq^{k-l+1}/y;q^{k-l+1};q,y)$
is analytic in $d$. Moreover, for $k\geq l$ we easily estimate
$|f_k(d)|\leq (-|d/y|,-|y|;q)_\iy$, so that
the convergence of the sum for $k\to\iy$ is uniform in $d$
on compact subsets since $|y|<1$. 

In order to obtain uniform convergence as $k\to -\iy$ we
use \thetag{\vglestimatektominusinfty} so that for $k\leq -l$
we can estimate $|f_k(d)|\leq |d^{l-k}(y/d;q)_{l-k}|
(-|d/y|,-|y|;q)_\iy$. On the other hand we see that 
by \thetag{\vglestimatektominusinfty} and 
\thetag{\vglthetaproduct} the coefficient for the function
$f_k(d)$ as $k\to\iy$ behaves like $C(-1)^kq^{\hf k(k+1)}$.
Hence the convergence for $k\to-\iy$ is also uniform in
$d$ on compact subsets. Hence both sides of the equality
of this proposition are analytic in $d$
and coincide on $|d|<1$, so it holds for all $d$.
\qed\enddemo

There are many ways of linking $w$ and $z$ such that the products
of the generating functions of 
Lemma \thmref{\lemgeneratingfunctions} simplify. Taking $w=x/z$
and $bx=y$ leads to the same result. The case $az=d/w$ also
gives an interesting result, which may be obtained as a limit
case of \cite{\KoelITSF, Prop.~2.2}.

In Proposition \thmref{\propfirststeporthorels} we 
replace $k$ by $k+p$, $l$ by $-r+p$ for $r,p\in\Z$, next
we take $y=q^{\al+1}$, $\al>-1$, and 
$aq^m=-cb$ to get
$$
\multline
\sum_{k=-\iy}^\iy {{q^{k(\al+1)}(q^{k+p+1};q)_\iy}
\over{(-cq^k;q)_\iy}}
 {}_1\vp_1\left( {{-cq^k}\atop{q^{k+p+1}}};
q,q^{\al+1} \right)
(q^{k+r+1};q)_\iy {}_1\vp_1\left( {{dq^{k+r-\al}}
\atop{q^{k+r+1}}};q,q^{\al+1} \right) \\
= \cases 0, &\text{for $p<r$,}\\
{\displaystyle{ q^{-p(\al+1)} d^{p-r}
{{(-cq^{\al-p+1}/d;q)_{p-r} (d,q^{p-r+1},q;q)_\iy}\over{
(-cq^{-p};q)_\iy}}}}, &\text{for $p\geq r$.}
\endcases
\endmultline
$$
Using \thetag{\vglthetaproduct} we see that this 
reduces to the second equation of
\thetag{\vglorthorelationscomplete} after specialising
$d=-cq^{\al-r}$, since $(q^{1+r-p};q)_{p-r}=\de_{r,p}$
for $p\geq r$. So we have proved 
\thetag{\vglorthorelationscomplete}.

%%%%%%%%%%%%%%%%%%%%%%%%%%%%%%%%%%%%%%%%%%%%%%%%%%%%%%%%%%%%%%%%%%%%
%%%%N E W   S E C T I O N%%%%%%%%%%%%%%%%%%%%%%%%%%%%%%%%%%%%%%%%%%%
%%%%%%%%%%%%%%%%%%%%%%%%%%%%%%%%%%%%%%%%%%%%%%%%%%%%%%%%%%%%%%%%%%%%
\subhead\newsection Proof of the dual orthogonality 
relations\endsubhead

In this section we prove the dual orthogonality relations
of Corollary \thmref{\eql}
by using a rigorous limit transition from the big $q$-Jacobi
polynomials to the $q$-Bessel functions under consideration,
which we will call now {\sl big $q$-Bessel functions}.
The orthogonality relations for the big $q$-Jacobi polynomials
tend to the dual orthogonality relations given in
Corollary \eql.
This then gives,
together with the analytic proofs of the previous section
a complete alternative analytic proof of
Theorem \thmref{\thmbasisappropriateLtwospace}.

For $\al>-1$, $c>0$, $k\in\Z$, $x\in\R$ put
$$
\align
\FSJ_{\al,k}^c(x;q):=&\;
{}_1\vp_1\left({x^{-1}\atop q^{\al+1}};q,-c^{-1}xq^{k+\al+2}\right)
\tag\eqname{\eqj}\\
=&\;(-c^{-1}q^{k+1};q)_\iy\
{}_2\vp_1\left({q^{\al+1}x,0\atop q^{\al+1}};
q,-c^{-1}q^{k+1}\right).
\tag\eqname{\eqc}
\endalign
$$
Here we used for the second identity the tranformation
$$
{}_1\vp_1\left({a\atop c};q,x\right)
=
(ax/c;q)_\iy\
{}_2\vp_1\left({c/a,0\atop c};q,ax/c\right),
$$
which is a limit case of Heine's
transformation formula \cite{\GaspR, (1.4.6)}.
Then formulas (\vgldefqLaguerrepol) and (\vgldefrestofbasis)
can be rewritten as
$$
L_n^{(\al)}(c^{-1}q^{k+1};q) =
{{(q^{\al+1};q)_n}\over{(q;q)_n}}\,
\FSJ_{\al,k}^c(q^n;q)\quad(n\in\Zp),
\tag\eqname{\eqa}
$$
respectively
$$
M_p^{(\al;c^{-1})}(c^{-1}q^{k+1};q) =
{(q^{\al+1};q)_\iy\over
(q,-c^{-1}q^{\al+1};q)_\iy}\,
\FSJ_{\al,k}^c(-cq^{p-\al};q)\quad(p\in\Z).
\tag\eqname{\eqb}
$$

The (dual) orthogonality relations in Corollary \eql\
can be rewritten by substitution of (\eqm), (\eqa) and (\eqj).
We obtain for $k,l\in\Z$:

$$
\multline
\sum_{n=0}^\iy (\FSJ_{\al,k}^c \FSJ_{\al,l}^c)(q^n;q)\,
{q^n\,(q^{n+1};q)_\iy\over(q^{n+\al+1};q)_\iy}\\
+\sum_{p=-\iy}^\iy (\FSJ_{\al,k}^c \FSJ_{\al,l}^c)
(-cq^{p-\al-1};q)\,
{cq^{p-\al-1}\,(-cq^{p-\al};q)_\iy\over(-cq^p;q)_\iy}\\
=\de_{k,l}q^{-k(\al+1)}\,{(q;q)_\iy^2\over(q^{\al+1};q)_\iy^2}\,
{(-cq^{-\al-1},-q^{\al+2}c^{-1},-q^{k+1}c^{-1};q)_\iy\over
(-c,-qc^{-1};q)_\iy}\,.
\endmultline
\tag\eqname{\eqe}
$$
We will show that the orthogonality relations (\eqe) 
can be rigorously
obtained as a limit case of the orthogonality relations for big
$q$-Jacobi polynomials.

Let $0<a<q^{-1}$, $b>-q^{-1}$, $c>0$.
{\sl Big $q$-Jacobi polynomials} are defined by
$$
P_k(x;a,b,-c;q):={}_3\vp_2\left({q^{-k},abq^{k+1},
x\atop aq,-cq};q,q\right)
\quad(x\in\R,\; k\in\Zp).
\tag\eqname{\eqf}
$$
Formally we have
$$
\lim_{r\to\iy}P_{r-k}(q^{\al+1}x;q^\al,b,-cq^{-r-1};q)=
{1\over (-c^{-1}q^{k+1};q)_\iy}\,\FSJ_k(x;\al,c;q).
\tag\eqname{\eqg}
$$
(Substitute (\eqf) and (\eqc) and take termwise limits
in the ${}_3\vp_2$.)$\;$
The orthogonality relations for big $q$-Jacobi polynomials 
are given by
$$
\multline
\int_{-cq}^{aq}(P_kP_l)(x;a,b,-c;q)\,
{(x/a,-x/c;q)_\iy\over(x,-bx/c;q)_\iy}\,d_qx\\
=\de_{k,l}\,M\,{1-abq\over1-abq^{2k+1}}\,
{(q,bq,-abq/c;q)_k\over(abq,aq,-cq;q)_k}\,(acq^2)^k\,q^{k(k-1)/2},
\endmultline
\tag\eqname{\eqh}
$$
where
$$
M:={(1-q)\,aq\,(q,-c/a,-aq/c,abq^2;q)_\iy\over
(aq,bq,-cq,-abq/c;q)_\iy}\,,
$$
see \cite{\GaspR, (7.3.12)--(7.3.14)}.
(Note the error on the \RHS\ of \cite{\GaspR, (7.3.13)}:
the factor $(-ac)^n$ must be replaced by $(-acq^2)^{-n}$.)

Formally the orthogonality relations (\eqe) can be obtained
as a limit case of the orthogonality relations (\eqh)
for big $q$-Jacobi polynomials.
In (\eqh) just replace $k$ by $r-k$, $l$~by $r-l$,
$c$ by $cq^{r-1}$,
and $a$ by $q^\al$, and let $r\to\iy$.
Because of these limit results we call the functions
$x\mapsto \FSJ_k(x;\al,c;q)$ {\sl big $q$-Besel functions}.

For $b=0$ we can transform the \RHS\ of
(\eqf) by the transformation formula \cite{\GaspR, 
Exercise 1.15(i)}:
$$
\Pt_k(x;a,0,-c;q):=(-c^{-1}q^{-k};q)_k\,P_k(x;a,0,-c;q)
=
{}_2\vp_1\left({q^{-k},aq/x\atop aq};q,-x/c\right).
\tag\eqname{\eqk}
$$
Formula (\eqg) for $b=0$ can also be obtained as a formal termwise
limit by substituting (\eqk) and (\eqj) and by taking termwise
limits in the ${}_2\vp_1$.

Fix $\al>-1$, $c>0$, and let $r\in\Zp$. From (\eqh) we 
get orthogonality
relations for the functions
$x\mapsto \Pt_{r-k}(q^{\al+1}x;q^\al,0,-cq^{-r-1};q)$\quad
($k=r,r-1,\ldots$):
$$
\multline
\sum_{n=0}^\iy
(\Pt_{r-k}\Pt_{r-l})(q^{n+\al+1};q^\al,0,-cq^{-r-1};q)\;
{q^n\,(q^{n+1},-c^{-1}q^{\al+n+r+2};q)_\iy\over
(q^{\al+n+1};q)_\iy}\\
+
\sum_{p=-r}^\iy
(\Pt_{r-k}\Pt_{r-l})(-cq^p;q^\al,0,-cq^{-r-1};q)\;
{cq^{p-\al-1}\,(-cq^{p-\al},q^{p+r+1};q)_\iy\over(-cq^p;q)_\iy}\\
=
\de_{k,l}\,{(q,-cq^{-r-\al-1},-q^{\al+r+2}/c;q)_\iy\over
(q^{\al+1},-cq^{-r};q)_\iy}\,
{(q,-c^{-1}q^{k+1};q)_{r-k}\over(q^{\al+1};q)_{r-k}}\,
q^{(\al+1)(r-k)}.
\endmultline
\tag\eqname{\eqi}
$$
Here $k,l=r,r-1,\ldots\;$.
Note that the orthogonality relations (\eqi) and (\eqe) 
have the same
structure. We want to show that we can take a rigorous 
limit for $r\to\iy$
of (\eqi) which yields (\eqe) preserving this structure.

\proclaim{Proposition \theoremname{\eqn}}
Fix $\al>-1$, $c>0$, $k\in\Z$. Then, for each $x\in\R$ we have the
pointwise limit
$$
\lim_{r\to\iy} \Pt_{r-k}(q^{\al+1}x;q^\al,0,-cq^{-r-1};q)=
\FSJ_{\al,k}^c(x;q).
\tag\eqname{\eqo}
$$
Furthermore, if $M>0$ then
$$
\multline
\left. \aligned
\bigl|\Pt_{r-k}(q^{\al+1}x;q^\al,0,-cq^{-r-1};q)\bigr|&\\
\bigl|\FSJ_{\al,k}^c(x;q)\bigr|&\endaligned\right\}
\le
{}_1\phi_1\left({-M^{-1}\atop q^{\al+1}};q,{-q^{\al+k+2}M
\over c}\right)\\
\hbox{for $|x|\le M$ and $r=k,k+1,\ldots\;$.}\endmultline
\tag\eqname{\eqp}
$$
\endproclaim

\demo{Proof}
Write
$$
\Pt_{r-k}(q^{\al+1}x;q^\al,0,-cq^{-r-1};q)=
{}_2\phi_1\left({q^{-r+k},x^{-1}\atop q^{\al+1}};q,-\,
{q^{\al+r+2}x\over c}
\right)=
\sum_{j=0}^\iy t_r(j,x)
$$
with
$$
t_r(j,x):={(q^{-r+k},x^{-1};q)_j\over(q^{\al+1},q;q)_j}\,
(-c^{-1}q^{\al+r+2}x)^j\quad
\hbox{(vanishing if $j>r-k$).}
$$
Also write
$$
\FSJ_{\al,k}^c(x;q)=
{}_1\phi_1\left({x^{-1}\atop q^{\al+1}};q,-c^{-1}
xq^{k+\al+2}\right)=
\sum_{j=0}^\iy t(j,x)
$$
with
$$
t(j,x):={(x^{-1};q)_j\over(q^{\al+1},q;q)_j}\,q^{\hf j(j-1)}\,
(c^{-1}xq^{k+\al+2})^j.
$$
Then
$$
\lim_{r\to\iy} t_r(j,x)=t(j,x)
\tag\eqname{\eqq}
$$
and
$$
|t_r(j,x)|\le T(j,M)\quad\hbox{for $r\ge k$ and $|x|\le M$,}
\tag\eqname{\eqr}
$$
where
$$
T(j,M):={(-M^{-1};q)_j\over(q^{\al+1},q;q)_j}\,q^{\hf j(j-1)}\,
(c^{-1}Mq^{k+\al+2})^j
$$
and
$$
\sum_{j=0}^\iy T(j,M)=
{}_1\phi_1\left({-M^{-1}\atop q^{\al+1}};q,-c^{-1}
Mq^{k+\al+2}\right)<\iy.
$$
Then the limit (\eqo) follows from (\eqq) and (\eqr) by dominated
convergence.

For the proof of (\eqr) we have used that, for 
$i=0,1,\ldots,r-k-1$ and
$|x|\le M$:
$$
\multline
\bigl|(1-q^{-r+k+i})(1-x^{-1}q^i)q^{r-k} x\bigr|=
\bigl|(q^i-q^{r-k})(x-q^i)\bigr|\\
\le
q^i(M+q^i)=(1+M^{-1}q^i)q^iM.\qquad\square
\endmultline
$$
\enddemo

For $\Pt_{r-k}(-cq^p;q^\al,0,-cq^{-r-1};q)$ as 
$-r\le p\le-k$, $r\to\iy$
we need a more refined estimate in order to be able to take limits
in (\eqi).

\proclaim{Proposition \theoremname{\eqs}}
Fix $\al>-1$, $c>0$, $k\in\Z$, Let $r\in\{k,k+1,\ldots\}$. 
Then, for $p\in\{-k,-k-1,\ldots\}$ we have
$$
\multline
\Pt_{r-k}(-cq^p;q^\al,0,-cq^{-r-1};q)=
(q;q)_{r-k}\,{(-cq^p;q)_{-k-p}\over(q^{\al+1},q;q)_{-k-p}}\,
(-c^{-1}q^{\al-p+1})^{-k-p}\\
\times
{}_2\phi_2\left({q^{-r-p},-cq^{-k}\atop q^{\al+1-k-p},q^{-k-p+1}};
q,-c^{-1}q^{\al-k-p+r+2}\right)
\endmultline
\tag\eqname{\eqt}
$$
and
$$
\multline
\FSJ_{\al,k}^c(-cq^{p-\al-1};\al,c;q)=
(q;q)_\iy\,{(-cq^p;q)_{-p-k}\over(q^{\al+1};q;q)_{-p-k}}\,
(-c^{-1}q^{-p+\al+1})^{-p-k}\\
\times{}_1\phi_2\left({-cq^{-k}\atop q^{-p-k+1},q^{\al+1-p-k}};q,
-c^{-1}q^{-2p-k+\al+2}\right)\\
=\lim_{r\to\iy}\Pt_{r-k}(-cq^p;q^\al,0,-cq^{-r-1};q).
\endmultline
\tag\eqname{\equ}
$$

Furthermore,
$$
\multline
\left. \aligned
\bigl|\Pt_{r-k}(-cq^p;q^\al,0,-cq^{-r-1};q)\bigr|&\\
\bigl|\FSJ_{\al,k}^c(-cq^{p-\al-1};\al,c;q)\bigr|&
\endaligned\right\}
\le
{(-cq^p;q)_{-k-p}\over(q^{\al+1},q;q)_{-k-p}}\,
(c^{-1}q^{\al-p+1})^{-k-p}\\
\times{}_1\phi_2\left({-cq^{-k}\atop q^{\al+1},q};q,c^{-1}
q^{k+\al+2}\right)\\
\hbox{for $r\in\{k,k+1,\ldots\}$ and $p\in\{-k,-k-1,\ldots\}$.}
\endmultline
\tag\eqname{\eqv}
$$
\endproclaim

\demo{Proof}
By Jackson's tranformation formula \cite{\GaspR, (1.5.4)} 
we have
$$
\align
{}_2\phi_1\left({a,b\atop c};q,z\right)&=
{(az;q)_\iy\over (z;q)_\iy}\
{}_2\phi_2\left({a,c/b\atop c,az};q,bz\right)\\
&={1\over(z;q)_\iy}\,
\sum_{j=0}^\iy{(a,c/b;q)_j\,(azq^j;q)_\iy\over(c,q;q)_j}\,
q^{\hf j(j-1)}(-bz)^j.
\endalign
$$
Put $z:=q^{-s+1}a^{-1}$, where $s\in\Zp$.
Then $(azq^j;q)_\iy=(q^{-s+j+1};q)_\iy$, so in the last 
sum the summation will start at $j=s$.
Replace the summation index $j$ by $j+s$ and write the 
resulting sum again as a ${}_2\phi_2$.
Then we obtain for $s\in \Zp$:
$$
{}_2\phi_1\left({a,b\atop c};q,a^{-1}q^{-s+1}\right)=
{(q;q)_\iy\,(c/b;q)_s\,b^s\over (qa^{-1};q)_\iy\,(c;q)_s}\
{}_2\phi_2\left({aq^s,cq^s/b\atop cq^s,q^{s+1}};q,qa^{-1}b\right).
$$
This yields (\eqt) (here $s=-k-p$).

The first identity in (\equ) can be derived in a similar 
way by starting with
$$
{}_1\phi_1\left({a\atop c};q,z\right)={1\over (c;q)_\iy}\,
\sum_{j=0}^\iy{(az/c;q)_j\,(q^jz;q)_\iy\over(q;q)_j}\,
q^{\hf j(j-1)}\,(-c)^j
$$
(this is (\eqaa) with the \RHS\ expanded)
and then putting $z:=q^{-s+1}$ ($s\in\Zp$).

Expand the \RHS\ of (\eqt). So
$$
\Pt_{r-k}(-cq^p;q^\al,0,-cq^{-r-1};q)=
\sum_{j=0}^\iy t_r(j;p)
$$
with
$$
\multline
t_r(j,p):=
(q;q)_{r-k}\,{(-cq^p;q)_{-k-p}\over(q^{\al+1},q;q)_{-k-p}}\,
(-c^{-1}q^{\al-p+1})^{-k-p}\\
\times{(q^{-r-p},-cq^{-k};q)_j\over(q^{\al+1-k-p},
q^{-k-p+1},q;q)_j}\,
q^{\hf j(j-1)}\,
(c^{-1}q^{\al-k-p+r+2})^j\quad
\hbox{(vanishing if $j>r+p$).}
\endmultline
$$
Then
$$
|t_r(j,p)|\le T(j,p)\quad
\hbox{for $r\in\{k,k+1,\ldots\}$ and $p\in\{-k,-k-1,\ldots\}$,}
\tag\eqname{\eqx}
$$
where
$$
T(j,p):={(-cq^p;q)_{-k-p}\over(q^{\al+1},q;q)_{-k-p}}\,
(c^{-1}q^{\al-p+1})^{-k-p}\,
{(-cq^{-k};q)_j\over(q^{\al+1},q,q;q)_j}\,q^{j(j-1)}\,
(c^{-1}q^{k+\al+2})^j
$$
and
$\sum_{j=0}^\iy T(j,p)$ equals the \RHS\ of (\eqv).

For the proof of (\eqx) we have used that, for $i=0,1,\ldots,r+p-1$,
$$
\bigl|(1-q^{-r-p+i})q^{r-p-k}\bigr|=
\bigl|(q^i-q^{r+p})q^{-2p-k}\bigr|\le q^i q^k.
$$

For fixed $p$, the limit formula in (\equ) follows from (\eqo),
but it follows also by taking a termwise limit for $r\to\iy$ on the
\RHS\ of (\eqt) and by using dominated convergence in view
of (\eqx).\qed
\enddemo

\LP
{\bf Proof of (\eqe).}\quad
Suppose $k\ge l$ (without loss of generality).
As $r\to\iy$, the \RHS\ of (\eqi) tends to the \RHS\ of (\eqe), 
and each term on the \LHS\ of (\eqi) tends to the corresponding 
term on the \LHS\ of (\eqe) (because of (\eqo)). We will show 
that the \LHS\ of (\eqi) tends rigorously to the \LHS\ of (\eqe) 
by splitting up the \LHS\ of (\eqi) as
$$
\sum_{n=0}^\iy+\sum_{p=-l+1}^\iy+\sum_{p=-k+1}^l+\sum_{p=-\iy}^{-k},
$$
and by using dominated convergence for each of the three 
infinite sums.

As for the first sum, by (\eqp) the $n^{\text{th}}$ term is 
bounded in
absolute value by $Cq^n$, where $C>0$ and independent of $n$.
Similarly, the $p^{\text{th}}$ term in the second sum is 
bounded in absolute value by $Cq^p$.

By (\eqv), the $p^{\text{th}}$ term in the fourth sum is 
bounded in absolute
value by
$$
C{(-cq^p;q)_{-k-p}\over(q^{\al+1},q;q)_{-k-p}}\,
{(-cq^p;q)_{-l-p}\over(q^{\al+1},q;q)_{-l-p}}\,
(c^{-1}q^{\al-p+1})^{-k-l-2p}
q^p\,{(-cq^{p-\al};q)_\iy\over(-cq^p;q)_\iy}\le
A\,q^{p^2}\,B^p,
$$
where $A,B,C>0$ and independent of $p$.\quad$\square$

\demo{Remark}
As we said in the 
Introduction the limit transition (\eqg) from big 
$q$-Jacobi polynomials
to big $q$-Bessel functions is inspired
by what happens at a quantum group level.

The regular representation of the quantum group $SU_q(2)$
on Podle\'s' quantum spheres naturally decomposes into
irreducibles, generated by suitable spherical functions
which turn out to depend on a single variable. Such functions
can thus be identified with ordinary big $q$-Jacobi polynomials,
see \cite{\NM}.

In the limit transition from $SU_q(2)$ to the Euclidean quantum
group $E_q(2)$, quantum spheres are replaced by quantum 
hyperboloids,
i.e., algebras in two generators $z$ and $\zb$ such that $z\zb=
q^2\zb z+1-q^2$. The corresponding regular representation 
decomposes 
into irreducibles, each of which has an $\iy$-dimensional basis
consisting of certain formal power series in $1-\zb z$, see
\cite{\BCG}.
This makes it
possible to identify such power series with ordinary functions,
and more precisely with $q$-Bessel functions of the form
$\FSJ_{\al,k}^c(x;q)$.

More precisely, the paper \cite{\BCG} gives in (4.9)
the basis elements of the
irreducible representation space as explicit formal power series
$J_r^{(q)}$. After the substitution
$\zb^j z^j=(1-\zb z;q^{-2})_j$ (see their Remark (4.10) (ii))
one can identify the series $J_r^{(q)}$ for ${\Cal E}=-q^{2k}$
with $\FSJ_{r,k}^{q^2/(q^2-1)^2}(1-\zb z;q^2)$ as defined by (\eqj).
\enddemo

%%%%%%%%%%%%%%%%%%%%%%%%%%%%%%%%%%%%%%%%%%%%%%%%%%%%%%%%%%%%%%%%%%%%
%%%%%R E F E R E N C E S%%%%%%%%%%%%%%%%%%%%%%%%%%%%%%%%%%%%%%%%%%%%
%%%%%%%%%%%%%%%%%%%%%%%%%%%%%%%%%%%%%%%%%%%%%%%%%%%%%%%%%%%%%%%%%%%%


\Refs

\ref\no \Akhi
\by N.I.~Akhiezer
\book The Classical Moment Problem and Some Related Questions in
Analysis
\publaddr Hafner 
\yr 1965
\endref 

\ref\no \Bere
\by J.M.~Berezanski\u\i
\book Expansions in Eigenfunctions of Selfadjoint Operators
\bookinfo Transl. Math. Monographs 17
\publaddr Amer. Math. Soc.
\yr 1968
\endref

\ref\no \Berg
\by C.~Berg
\paper On some indeterminate moment problems for measures
on a geometric progression
\paperinfo preprint
\yr 1997
\endref

\ref\no \BCG
\by  F.~Bonechi, N.~Ciccoli, R.~Giachetti, E.~Sorace and M.~Tarlini 
\paper Free $q$-Schr\"odinger equation from quantum 
homogeneous spaces of the 2-dim Euclidean quantum group
\jour Comm. Math. Phys.
\vol 175
\yr 1996
\pages 161--176
\endref

\ref\no \ChenIM
\by Y.~Cheng, M.E.H.~Ismail and K.A.~Muttalib
\paper Asymptotics of basic Bessel functions and $q$-Laguerre
polynomials
\jour J. Comp. Appl. Math.
\vol 54
\yr 1994
\pages 263--272
\endref

\ref\no \DunfS
\by N.~Dunford and J.T.~Schwartz
\book Linear Operators II: Spectral Theory
\publaddr Interscience
\yr 1963
\endref

\ref\no \GaspR
\by G.~Gasper and M.~Rahman
\book Basic Hypergeometric Series
\publaddr Cambridge Univ. Press
\yr 1990
\endref

\ref\no \Isma
\by M.E.H.~Ismail
\paper The zeros of basic Bessel functions, the functions
$J_{\nu+ax}(x)$, and associated orthogonal polynomials
\jour J. Math. Anal. Appl.
\vol 86
\yr 1982
\pages 1--19
\endref

\ref\no \IsmaMS
\by M.E.H.~Ismail, D.R.~Masson and S.K.~Suslov
\paper The $q$-Bessel function on a $q$-quadratic grid
\paperinfo preprint
\yr 1996
\endref

\ref\no \IsmaR
\by M.E.H.~Ismail and M.~Rahman
\paper The $q$-Laguerre polynomials and related moment problems
\jour J. Math. Anal. Appl.
\vol 218 \yr 1998 \pages 155--174
\endref

\ref\no \Kake
\by T.~Kakehi
\paper Eigenfunction expansion associated with the Casimir
operator on the quantum group $SU_q(1,1)$
\jour Duke Math. J.
\vol 80
\yr 1995
\pages 535--573
\endref

\ref\no \KakeMU
\by T.~Kakehi, T.~Masuda and K.~Ueno
\paper Spectral analysis of a $q$-difference operator which 
arises from the quantum $SU(1,1)$ group
\jour J. Operator Theory
\vol 33
\yr 1995
\pages 159--196
\endref

\ref\no \KoelITSF
\by H.T.~Koelink
\paper A basic analogue of Graf's addition formula and
related formulas
\jour Integral Transforms and Special Functions
\vol 1
\yr 1993
\pages 165--182
\endref

\ref\no \KoelDMJ
\bysame % H.T.~Koelink
\paper The quantum group of plane motions and the Hahn-Exton 
$q$-Bessel function
\jour Duke Math. J.
\vol 76
\yr 1994
\pages 483--508
\endref

\ref\no \KoelIM
\bysame % H.T.~Koelink
\paper The quantum group of plane motions and basic Bessel functions
\jour Indag. Mathem. N.S.
\vol 6
\yr 1995
\pages 197--211
\endref

\ref\no \KoelV
\by H.T.~Koelink and J.~Verding
\paper Spectral analysis and the Haar functional on the 
quantum $SU(2)$ group
\jour Comm. Math. Phys.
\vol 177
\yr 1996
\pages 399--415
\endref

\ref\no \Koor
\by T.H.~Koornwinder
\paper Jacobi functions and analysis on noncompact
semisimple Lie groups
\inbook Special Functions: Group Theoretical Aspects 
and Applications
\eds R.A.~Askey, T.H.~Koornwinder, W.~Schempp
\pages 1--85
\yr 1984
\publaddr Reidel
\endref

\ref\no \KoorS
\by T.H.~Koornwinder and R.F.~Swarttouw
\paper On $q$-analogues of the Fourier and Hankel transforms
\jour Trans. Amer. Math. Soc.
\vol 333
\yr 1992
\pages 445--461
\endref

\ref\no \MassR
\by D.R.~Masson and J.~Repka
\paper Spectral theory of Jacobi matrices in $\ell^2(\Z)$ 
and the ${\frak {su}}(1,1)$ Lie algebra
\jour SIAM J. Math. Anal. \vol 22 \yr 1991
\pages 1131--1146
\endref

\ref\no \MMNNSU
\by T.~Masuda, K.~Mimachi, Y.~Nakagami, M.~Noumi, Y.~Saburi 
and K.~Ueno
\paper Unitary representations of the quantum group $SU_q(1,1)$: 
II---Matrix elements of unitary representations and the basic 
hypergeometric functions
\jour Lett. Math. Phys.
\vol 19
\yr 1990
\pages 195--204
\endref

\ref\no \Moak 
\by D.S.~Moak
\paper The $q$-analogue of the Laguerre polynomials
\jour J. Math. Anal. Appl. \vol 81 \yr 1981 \pages 20--47
\endref

\ref\no \NM
\by M.~Noumi and K.~Mimachi
\paper Quantum 2-spheres and big $q$-Jacobi polynomials
\jour Comm. Math. Phys.
\vol 128
\yr 1990
\pages 521--531
\endref

\ref\no \VK
\by  L.L.~Vaksman and L.I.~Korogodski{\u\i}
\paper An algebra of bounded functions on the quantum group of 
the motions of the plane and $q$-analogues of Bessel functions
\jour Soviet Math. Dokl.
\vol 9
\yr 1989
\pages 173-177
\endref

\endRefs
\enddocument